
\magnification=1200
\input amstex
\documentstyle{amsppt}
\NoBlackBoxes
\NoRunningHeads
\TagsOnRight

\define\C{{\Bbb C}}
\define\Z{{\Bbb Z}}
\define\De{{\Delta}}
\define\ra{{\rightarrow}}
\define\a{{\alpha}}

\topmatter
\title
On quantum group $GL_{p,q}(2)$\endtitle

\author
{Tanya  Khovanova}
\endauthor

\address{\it Present address:} Tanya Khovanova, 
Department of
Mathematics, Princeton University, NJ 08544
\endaddress

\email{tanya\@math.princeton.edu} \endemail 

\abstract
In the Hopf algebra $Gl_{p,q}(2)$ the determinant is central
iff $p=q$. In this case we put determinant to be equal to $1$
to get $SL_q(2)$. In this paper we consider the case when 
$p/q$ is a root of unity; and, consequently, a power
of the determinant is central.
\endabstract

\date January 28, 1997
\enddate

\endtopmatter

\document
\subhead
{Introduction}
\endsubhead
\medskip
In the paper [B-Kh] the Universal enveloping algebras for
metaplectic quantum groups of $SL(2)$-type were
constructed. As we look back at the history
of $SL_q(2)$ we can notice that on the levels
of formulas the algebra of functions on $SL_q(2)$ was always
presented the same way ([M], [R-T-F]), while the Universal enveloping
algebra had different presentations. This fact encouraged
me and Joseph Bernstein [B-Kh] to consider 
the enveloping algebra as the secondary object in
comparison with the Hopf algebra of functions.
As a byproduct of our axiomatic search we
constructed the metaplectic quantum groups of $SL(2)$-type.
This construction was carried out on the level
of universal algebras. 
\medskip
On my way to find nice formulas for
algebra of functions on metaplectic quantum group, I found
some nice formulas for other algebras. These formulas
are presented in this paper.
\medskip
I am thankful to Joseph Bernstein for his interest
in this subject and helpful discussions.
\bigskip
\subhead
{1. Algebra of regular functions on $GL_{p,q}(2)$}
\endsubhead
\medskip
{\bf 1.1.} For every $p,q \in \C^*$ we consider the bialgebra $Mat_{p,q}$
which is 
generated by four noncommuting elements $(a,b,c,d)$, satisfying the
following
relations (see, [M],  [O-W] and references there):
$$\aligned 
ab & = p^{-1}ba \\
cd & = p^{-1}dc \\
bc & = q^{-1}pcb
\endaligned 
\qquad
\aligned
ac &=q^{-1}ca \cr
bd &=q^{-1}db \cr
ad-da &=(p^{-1}-q)bc.
\endaligned \tag *$$
\medskip
Introduce matrices
$$ \align 
Y &= \pmatrix a & b \cr c & d \endpmatrix \in Mat\ (2,Mat_{p,q}) \cr \cr
P &= \pmatrix
0 & -1 \cr p^{-1} & 0 \endpmatrix \qquad Q = \pmatrix 0 & -1 \cr q^{-1}
& 0 \endpmatrix.\endalign $$
Then we can rewrite the relations $(*)$ in a more compact form:
$$ \eqalign {
Y P Y^t &= DP \cr
Y^t Q Y &= DQ, \cr
}$$
where $D$ is an element in $Mat_{p,q}$ and has a meaning
of quantum determinant.
\medskip
{\it Comment.} The matrix $P$ defines a quantum plane as
an algebra generated by two generators $x$, $y$ and the
relation $xy=p^{-1}yx$. Analogously the matrix $Q$ defines
a quantum plane. We consider $Mat_{p,q}$ as an algebra of operators
preserving quantum plane defined by $P$ and dual quantum plane
defined by $Q$.
\medskip
{\bf 1.2.} 
The comultiplication in the algebra $Mat_{p,q}$ is defined as follows:
$$\aligned 
\De a &= a \otimes a +  b \otimes c    \cr
\De b &= a \otimes b +  b \otimes d    \cr
\De c &= c \otimes a +  d \otimes c    \cr
\De d &= c \otimes b +  d \otimes d \ .
\endaligned \tag **$$
Using the natural imbeddings $i',\ i'' : Mat_{p,q} \ra Mat_{p,q}
\otimes Mat_{p,q},\ \bigl (
i' (x) = x \otimes 1,\ i''(x) = 1 \otimes x \bigr )$, we can rewrite
comultiplication formulae $(**)$ as follows:
$$\De(Y) = i'(Y) \cdot i''(Y)\ ,$$
which is an equality in $Mat (2,Mat_{p,q} \otimes Mat_{p,q})$.
\medskip
{\bf 1.3.}
This algebra has a multiplicative quantum determinant 
$D=det_{p,q}(Y)$ (see [M], [O-W], [K]):
$$D= da -pcb=da -qbc=ad -p^{-1}bc=ad -q^{-1}cb.$$
Multiplicativity means that $\De D=D \otimes D$, or,
equivalently, $det_{p,q}(Y_1Y_2) = det_{p,q}(Y_1)\cdot
det_{p,q}(Y_2)$
whenever the entries of $Y_1$ commute with the entries of $Y_2$.
\medskip
Localizing by $D^{-1}$ we will get the algebra of functions on
$GL_{p,q}(2)$. This localization can be described easily due
to the fact that $D$ is normalizing (see [M]):
$$\align 
Da & =aD\qquad Db=p^{-1}qbD \cr
Dc & =pq^{-1}cD \qquad Dd=dD.
\endalign$$
\medskip
{\bf 1.4.}
The bialgebra $Fun(GL_{p,q})$ is the Hopf algebra;
an antipode $S$ is defined by:
$S(Y)= Y^{-1}$. Specifically:
$$\align
S(a)& =dD^{-1} \qquad S(b)=-pbD^{-1} \cr
S(c)&= -p^{-1}cD^{-1}\qquad
S(d)=aD^{-1}.
\endalign$$
In a more compact form:
$$Y^{-1}=
 S(Y)= PY^tP^{-1}D^{-1}=D^{-1}Q^{-1}Y^tQ.$$
\medskip
{\bf 1.5.}
If $p=q$ then $D$ is central; and we can take quotient of
$Fun(GL_{p,q}(2))$ by the Hopf ideal generated by $D-1$. We would get 
standard $Fun(SL_q(2))$ (see [Kas]).
\medskip
{\bf 1.6.}
Suppose now 
that $p^{-1}q$ equals $\xi$, where $\xi$ is the n-th
root of unity: $\xi^n=1$. In this case $D^n$ is central. 
An ideal generated by $D^n-1$ is a Hopf ideal, so
we can consider a quotient algebra, which we would denote by
$Fun(SL_{q,\xi}(2))$:
$$Fun(SL_{q,\xi}(2)) = Fun(GL_{p,q}(2))/(D^n-1).$$
\bigskip
\subhead{2. Axiomatic approach}\endsubhead
\medskip
{\bf 2.1.} Analyzing the Hopf algebra $A = Fun(SL_{q,\xi}(2))$ 
we note that it has the following
important property:
\medskip
Let $I \subset A$ be a two-sided ideal generated by $b $ and $c$. 
Then $I$ is a
Hopf ideal in $A$, i.e. $\De I \in A \otimes I + I \otimes A$ and
$S(I) \subset I$.
The quotient Hopf algebra $A/I$ is
isomorphic to the algebra of functions on the algebraic group $\C^* 
\otimes \Z_n$:
$$\C [a,d,D,D^{-1}] \big  / (ad-D,D^n-1) $$
$$\De a = a \otimes a \ \ \ \De d = d \otimes d\ \ \ \De D =
D \otimes D$$
$$S(a)=dD^{-1}\ \ \ S(d)=aD^{-1}\ \ \ S(D)=D^{-1}.$$
Or, equivalently:
$$\C [a,a^{-1},D,D^{-1}
] \big / (D^n-1)$$
$$\De a = a \otimes a \qquad \De D = D \otimes D$$
$$S(a)=a^{-1}\ \ \ S(D)=D^{-1}.$$
Informally, this means that our quantum group $A=Fun(SL_{q,\xi}(2))$ 
contains
the group $\C^* \otimes \Z_n$ as a subgroup.
\bigskip
\subhead{3. Dual picture}\endsubhead
\medskip
{\bf 3.1.}
Let us denote by $T^2$ the two-dimensional torus. Denote by $t$ a
point of $T^2$. To each point $t=(x_1,x_2)$ we can correspond a generator 
$\hat t$ in the group algebra of torus. This correspondence is
multiplicative --- $\widehat {t_1t_2}=\hat t_1\hat t_2$.
\medskip
We can describe (we use notations similar to [B-Kh]) a universal 
enveloping algebra $U_{p,q}(2)$ of $GL_{p,q}(2)$ as a Hopf algebra 
generated by $\hat t\ (t \in T^2)$ and two elements $E, F$,
satisfying the relations:
$$\aligned
\hat t E \hat t^{-1} & = \a(t) E \cr
\hat t F \hat t^{-1} & = (-\a)(t) F 
\endaligned \tag 1$$
$$ \aligned
\De \hat t & = \hat t \otimes \hat t \cr
\De E & = E \otimes 1 + Q_1 \otimes E \cr
\De F & = F \otimes Q_2^{-1} + 1 \otimes F 
\endaligned \tag 2$$
$$
[E,F] =  \frac {Q_1 - Q_2^{-1}} {q - p^{-1}}.
\leqno (3)$$
Here $Q_1, Q_2$ are generators corresponding to points of
the torus with coordinates $(q,p^{-1})$ and $(p,q^{-1})$
respectively; and $\a$ denotes the weight $(1,-1)$ on torus:
if $t=(x,y)$, then $\a(t)=xy^{-1}$.
\medskip
{\it Comments.} 1. The commutator $[E,F]$ is the element $X$ of an algebra,
generated by torus, satisfying the equation: $\De X=
Q_1 \otimes X + X \otimes Q_2^{-1}$. 
In denominator (3) we can choose any constant.
We chose our constant so that the evaluation of the right hand side element
on the weight $(1,0)$ would be equal to $1$.
\medskip
2. In this case and in the following cases, the antipode is
uniquely defined and could be easily recovered:
$$S(\hat t) = \hat t^{-1}\ \ \ S(E) = - Q^{-1}_1E \ \ \ S(F)=-Q_2F.$$
{\bf 3.2.}
To get the universal algebra $U_q(2)$ of $SL_q(2)$ we have
to put $p=q$ and to take the group subalgebra of $T^1 \subset T^2$ generated
by elements $\hat h$ corresponding to
 $(h,h^{-1})\in T^2$. Denote by $K$ an element
corresponding to $(q,q^{-1})$, then we would have the
following relations:
$$\leqalignno {
\hat h E \hat h^{-1} & = h^2 E &(1)\cr
\hat h F \hat h^{-1} & = h^{-2} F \cr
}$$
$$ \leqalignno {
\De \hat h & = \hat h \otimes \hat h &(2)\cr
\De E & = E \otimes 1 + K \otimes E \cr
\De F & = F \otimes K^{-1} + 1 \otimes F \cr
}$$
$$
[E,F] =  \frac {K - K^{-1}} {q - q^{-1}}.
\leqno (3)$$
\medskip
{\bf 3.3.}
Consider the case $SL_{q,\xi}(2)$: namely, $p^{-1}q=\xi$,
where $\xi^n=1$. We can describe the universal enveloping
algebra $U_{q,\xi}(2)$ as a subalgebra in $U_{p,q}(2)$, generated
by elements $E,\ F,\ \hat h=(h, h^{-1}) \in T^2$ for $h\in T$ and 
$W=(1,\xi)\in T^2$ ($W^n=1$).
$$\leqalignno {
WEW^{-1}= \xi^{-1}E \qquad \hat h E \hat h^{-1} & = h^2 E &(1)\cr
WFW^{-1}= \xi F \qquad \hat h F \hat h^{-1} & = h^{-2} F \cr
}$$
$$ \leqalignno {
\De W & = W\otimes W\cr
\De \hat h & = \hat h \otimes \hat h &(2)\cr
\De E & = E \otimes 1 + W\hat q \otimes E \cr
\De F & = F \otimes W \hat \xi \hat q^{-1} + 1 \otimes F \cr
}$$
$$
[E,F] =  \frac {\hat q - \hat \xi \hat q^{-1}} {q -\xi q^{-1}}W.
\leqno (3)$$
\bigskip
\bigskip
\noindent {\bf References}
\medskip
[B-Kh] J.Bernstein and T.Khovanova, \it On quantum group 
$SL_q(2)$, \rm Comm. Math. Phys., \bf 117 \rm (1996),
691-708. {\it hep-th/9412056}.
\medskip
[K] B.A.Kupershmidt, \it Classification of quantum group
structures on the group $GL(2)$, \rm J. Phys. A., \bf 27 \rm (1994),
L47-L51.
\medskip
[Kas] Ch.Kassel, \it Quantum groups, \rm Springer-Verlag. 1995.
\medskip
[M] Yu.I.Manin,  \it Topics in Noncommutative Geometry,
\rm Princeton University Press, 1991.
\medskip
[O-W] O.Ogievetsky, J.Wess, \it Relations between
$GL_{p,q}(2)$'s, \rm Z. Phys. C, \bf 50 \rm (1991), 123-131.
\medskip
[R-T-F] N.Yu.Reshetikhin, L.A.Takhtadzhyan and L.D.Faddeev,
\it Quantization of Lie groups and Lie algebras, \rm 
Leningrad Math. J. \bf 1 \rm (1990), 193-225.
\enddocument